\documentclass{elsarticle}

\usepackage{amsmath}
\usepackage{amsfonts}
\usepackage{amssymb}

\newcommand{\abs}[1]{\ensuremath{\left| #1 \right| }}
\newcommand{\qPs}{$q$--Pochhammer symbol}
\newcommand{\bhs}{basic hypergeometric series}

\newcommand{\lhs}{left hand side}

\newcommand{\od}{one dimensional}

\newcommand{\RSi}{Rogers--Selberg identity}

\newcommand{\hgr}{hyperoctahedral group}
\newcommand{\Js}{Jackson sum}

\newcommand{\tns}{\ensuremath{_{10}\varphi_9}}

\newcommand{\sss}{\ensuremath{_6\psi_6}}

\newdefinition{defn}{Definition}
\newdefinition{rmk}{Remark}
\newproof{pf}{Proof}

\begin{document}

\begin{frontmatter}

\journal{Adv and App in Discete Math}
\date{February 21, 2010}

\title{Multilateral basic hypergeometric summation identities \\ and hyperoctahedral group symmetries}

\author{Hasan Coskun}
\ead{hasan\_coskun@tamu-commerce.edu}
\ead[url]{http://faculty.tamu-commerce.edu/hcoskun}
\address{Department of Mathematics, Texas A\&M
  University--Commerce, Binnion Hall, Room 314, Commerce, TX 75429}  

\begin{abstract}
We give new proofs for certain bilateral basic hypergeometric summation formulas using the symmetries of the corresponding series. In particular, we present a proof for Bailey's $_3\psi_3$ summation formula as an application. 
We also prove a multiple series analogue of this identity by considering 
 hyperoctahedral group symmetries of higher ranks.
\end{abstract}

\end{frontmatter}

\section{Introduction}
Let $(a;q)_\alpha$ denote the \qPs\ which is formally defined by 
\begin{equation}
\label{qPochSymbol} 
(a)_\alpha = (a;q)_\alpha :=\dfrac{(a;q)_\infty}{(aq^\alpha;q)_\infty}
\end{equation} 
where the parameters $a, q, \alpha\in\mathbb{C}$, and $(a;q)_\infty$ denotes the infinite product $(a;q)_\infty:=\prod_{i=0}^{\infty} (1-aq^i)$. Note here that when $\alpha=k$ is a positive integer, then the \qPs\ reduces to $(a)_k = \prod_{i=1}^{k} (1-aq^{i-1}) $. We often use the shorthand notation $(a_1, a_2, \ldots, a_r)_\alpha$ for the product $\prod_{i=1}^{r} (a_i)_\alpha$.

The series $\sum_{k=0}^\infty c_k$, where the ratio $c_{k+1}/c_{k}$ is
a rational function of $q^k$, is called a basic hypergeometric
series~\cite{Andrews0}. Using the \qPs~(\ref{qPochSymbol}), the general basic
hypergeometric series with $r$ numerator parameters and $s$
denominator parameters is defined by 
\begin{equation}
  \label{eq:BHS}
_r\varphi_s\left[
\begin{matrix} 
a_1,a_2,\ldots,a_r \\
b_1, b_2, \ldots, b_s
\end{matrix}
;x , q\right] 
:=\sum_{k=0}^\infty \dfrac{(a_1, a_2, \ldots , a_r)_k }{(q, b_1,
  b_2, \ldots, b_s)_k }  \, x^k \left((-1)^k q^{\binom{k}{2}}\right)^{1+s-r}
\end{equation}
where we assume that none of the denominator factors vanish.

Note that if one of the numerator parameters is of
the form $q^{-n}$, for some non--negative integer $n$ and $q\neq 0$,
the series terminates from above since $(q^{-n})_k=0$ when $k>n$. The
denominator factor $(q)_k$ terminates the 
series from below, that is the factor $1/(q)_n$ causes the sum to
vanish when $n<0$. 

In general, when dealing with non-terminating series it is assumed for convergence that
$\abs{q}<1$. In that case, the series $_{r+1}\varphi_r$ converges
absolutely for $\abs{x}<1$.

When $r=s+1$, the basic hypergeometric series~(\ref{eq:BHS})
is called well--poised if the parameters satisfy the relation 
\[qa_1=a_2 b_1 =a_3 b_2 =\ldots = a_{s+1} b_s, \]
and very well--poised if, in addition, $a_2 = q\sqrt{a_1}$ and
$a_3 = -q\sqrt{a_1}$. An $_{r+1}\varphi_r$ series is called
$k$-balanced if $b_1\ldots b_r = q^k a_1\ldots a_{r+1}$, and
$x=q$.  

There are numerous classical one-dimensional results, summation and
transformation formulas for \bhs. One of the most general summation
formulas, for example, is the ($q$--Dougall or) \Js\
\begin{multline}
\label{OneDimJs}
_8\varphi_7 \left[
\begin{matrix} 
a, q a^{1/2}, -q a^{1/2}, b, c ,d, e, q^{-n} \\
a^{1/2}, -a^{1/2}, aq/b, aq/c, aq/d, aq/e, aq^{n+1}
\end{matrix}
;q , q \right] \\ = \dfrac{(aq, aq/bc, aq/bd,
aq/cd)_n} {(aq/b, aq/c, a q/d, a q/bcd)_n } \hspace*{0.3in}
\end{multline}
where $qa^2=bcdeq^{-n}$. An important general transformation
formula is Bailey's \tns\ transformation
\begin{multline}
\label{OneDimTns}
_{10}\varphi_9 \left[
\begin{matrix} 
a, q a^{1/2}, -q a^{1/2}, b, c ,d, e, f, \lambda a q^{n+1}/ef, q^{-n} \\
a^{1/2}, -a^{1/2}, aq/b, aq/c, aq/d, aq/e, aq/f, ef q^{-n}/\lambda, a
q^{n+1}  
\end{matrix}
;q , q \right] \\ = \dfrac{(aq, aq/ef, \lambda q/e,
\lambda q/f)_\infty} {(aq/e, aq/f, \lambda q, \lambda q/ef)_\infty } \\
\cdot \, {_{10}\varphi_9} \left[
\begin{matrix} 
\lambda, q \lambda^{1/2}, -q \lambda^{1/2}, \lambda b/a, \lambda c/a,
\lambda d/a, e, f, \lambda a q^{n+1}/ef, q^{-n} \\
\lambda^{1/2}, -\lambda^{1/2}, aq/b, aq/c, aq/d, \lambda q/e, \lambda
q/f, ef q^{-n}/a, \lambda q^{n+1}  
\end{matrix}
;q , \dfrac{aq}{ef} \right]
\end{multline}
where $\lambda = qa^2/bcd$.

The \bhs~(\ref{eq:BHS}) is Heine's generalization of the
hypergeometric series 
\begin{equation}
  \label{eq:HS}
_r F_s\left[
\begin{matrix} 
a_1,a_2,\ldots,a_r \\
b_1, b_2, \ldots, b_s
\end{matrix}
;x \right] 
=\sum_{n=0}^\infty \dfrac{ \{a_1\}_n \{a_2\}_n \ldots \{a_r\}_n }{n!
  \, \{b_1\}_n \{b_2\}_n \ldots \{b_s\}_n } x^n
\end{equation}
where $\{a \}_n$ denotes the shifted factorial (or Pochhammer symbol)
defined by 
\begin{equation}
 \label{def:Ps}
\{a\}_0:=1, \quad \{a\}_n:=a(a+1)\cdots (a+n-1) \; \mathrm{for} \;
n\in\mathbb{Z}_{>}.
\end{equation}
The \bhs~(\ref{eq:BHS})
reduces to~(\ref{eq:HS}) if we replace parameters  $a_i$ and $b_i$ by
$q^{a_i}$ and $q^{b_i}$ in~(\ref{eq:BHS}) respectively, and let
$q\rightarrow 1$.

The \bhs\ are further generalized in the literature in several
directions. Bilateral 
\bhs\ is a generalization where the index of summation is no 
longer restricted to non--negative integers, but it runs over all
integers. The most general result of this type is Bailey's \sss\
summation formula which can be written as
\begin{multline}
_6\psi_6 \left[
\begin{matrix} 
q a^{1/2}, -q a^{1/2}, b, c ,d, e \\
a^{1/2}, -a^{1/2}, aq/b, aq/c, aq/d, aq/e
\end{matrix}
;q , \dfrac{qa^2}{bcde} \right] \\
= \dfrac{(aq, aq/bc, aq/bd, aq/be, 
aq/cd, aq/ce, aq/de, q, q/a)_\infty} {(aq/b, aq/c, a q/d, aq/e, q/b,
q/c, q/d, q/e, qa^2/bcde)_\infty }
\end{multline}
provided that $\abs{qa^2/bcde}<1$, where 
\begin{multline}
  \label{eq:bilateralBHS}
_r\psi_s\left[
\begin{matrix} 
a_1,a_2,\ldots,a_r \\
b_1, b_2, \ldots, b_s
\end{matrix}
;x , q\right] \\
=\sum_{n=-\infty}^\infty \dfrac{(a_1)_n (a_2)_n \ldots (a_r)_n }{ (b_1)_n
  (b_2)_n \ldots (b_s)_n } (-1)^{(s-r)n}
  q^{(s-r)\binom{n}{2}} x^n \hspace*{0.5in}
\end{multline}
There are other important summation formulas such as Ramanujan's
$_1\psi_1$ sum, and useful transformation formulas for bilateral
series as well. The former identity, for example, may be written as
\begin{equation}
  \label{1psi1}
_1\psi_1(a,b;x , q) \\
=\sum_{n=-\infty}^\infty \dfrac{(a)_n }{ (b)_n } \, x^n = 
\dfrac{(q, b/a, ax, q/ax)_\infty} {(b, q/a, x, b/ax)_\infty }
\end{equation}

Appell and Lauricella series~\cite{Slater1} are other generalizations
of the \bhs, where the number of variables in the argument is extended to
several variables and thus multisums are considered.

Macdonald~\cite{Macdonald3} generalized the \bhs\ to a multiple
\bhs\ of symmetric function argument
where the argument is replaced by a Schur function, which is the ratio
of two determinants, and the index of summation runs over
partitions.   
Several summation and transformation results have been obtained at this
generality as well. But all such results generalized lower
level identities which satisfy only the balancedness condition. 

The study of elliptic (modular) analogue of \bhs\ was started
by Frenkel and Turaev~\cite{FrenkelT1} who, 
defining an elliptic analogue of the \qPs, proved an elliptic
analogue of the \tns\ transformation given above. Several other
authors contributed by proving various \od\ and a few multisum
identities including~\cite{Warnaar2}
and~\cite{RosengrenS1}.

The elliptic generalization of the \qPs\ is given by means of the
normalized elliptic function
\begin{equation}
  \theta(x) = \theta(x;p) := (x; p)_\infty (p/x; p)_\infty  
\end{equation}
where $\abs{p}<1$. The elliptic \qPs\ is then given by 
\begin{equation}
  (a; q,p)_n = \prod_{k=0}^{n-1} \theta(aq^k) 
\end{equation}
for $n>0$. The definition is extended to negative $n$ by the relation
  $(a; q,p)_n = 1/ (aq^{n}; q, p)_{-n} $ analogous to the standard 
  \qPs. When $n=0$, we have $(a; q,p)_n=1$. Note also that when $p=0$
  this reduces to standard definition of the \qPs.

The definition of a balanced,
very--well--poised elliptic \bhs\ now may be written~\cite{Warnaar2} as
\begin{equation}
{\!\!\!}_{r+1}\omega_r (a_1; a_4, \ldots, a_{r+1}; q, p) = \sum_{k=0}^\infty
\dfrac{\theta(a_1q^{2k}) } {\theta(a_1) } \dfrac{(a_1, a_4, \ldots, a_{r+1}; q,
  p)_k q^k } {(q, a_1q/a_4, \ldots, a_1q/a_{r+1}; q, p)_k }
\end{equation}
where $(a_4\ldots a_{r+1})^2 = a_1^{r-3} q^{r-5}$. By defining the
partition generalization of the elliptic \qPs\ in the form
\begin{equation}
\label{ellipticqtPs}
  (a)_\lambda= (a; q,p, t)_\lambda := \prod_{k=0}^{n-1} (at^{1-i}; q,p)_{\lambda_i}
\end{equation}
the definition of elliptic \bhs\ is generalized to various root
systems of rank $n$. The following shorthand notation will also be used.
\begin{equation}
  (a_1, \ldots, a_k)_\lambda = (a_1, \ldots, a_k; q, p, t)_\lambda :=
  (a_1)_\lambda \ldots (a_k)_\lambda .
\end{equation}

In~\cite{CoskunG1} we proved a multiple elliptic analogue of the classical \Js\ and other important results including a multiple analogue of Bailey's \tns\ transformation formula. The multiple elliptic Jackson sum may be written in the form
\begin{multline}
\label{WJackson}
W_{\lambda}(z;q,p,t,at^{-2n},bt^{-n})\\
=\dfrac{(s)_{\lambda}(as^{-1}t^{-n-1})_{\lambda}}
{(qbs^{-1}t^{-1})_{\lambda}(qbt^ns/a)_{\lambda}}
\cdot \prod_{1\leq i < j\leq n}
\left\{\dfrac{(t^{j-i+1})_{\lambda_i-\lambda_j}
(qbt^{-i-j+1})_{\lambda_i+\lambda_j}}
{(t^{j-i})_{\lambda_i-\lambda_j} (qbt^{-i-j})_{\lambda_i +
    \lambda_j}} \right\}  \\
\cdot\sum_{\mu \subseteq \lambda}
  \dfrac{(bs^{-1}t^{-n})_{\mu} (qbt^n/a)_\mu}{(qt^{n-1})_\mu
(as^{-1}t^{-n-1})_\mu} \cdot \prod_{i=1}^n \left\{
  \dfrac{(1-bs^{-1}t^{1-2i}q^{2\mu_i})}{(1-bs^{-1}t^{1-2i})}(qt^{2i-2})^{\mu_i}
\right\} \\ \cdot \prod_{1\leq i < j \leq n}
\left\{ \dfrac{(t^{j-i})_{\mu_i -\mu_j} (qt^{j-i})_{\mu_i
      -\mu_j}}{(qt^{j-i-1})_{\mu_i -\mu_j}(t^{j-i+1})_{\mu_i
      -\mu_j}} \dfrac{(bs^{-1}qt^{-i-j})_{\mu_i+\mu_j}
    (bs^{-1}t^{-i-j+2})_{\mu_i+\mu_j}} {(bs^{-1}t^{-i-j+1})_{\mu_i+\mu_j}
(qbs^{-1}t^{-i-j+1})_{\mu_i+\mu_j}} \right\} \\
\cdot W_\mu(q^\lambda t^{\delta(n)};q,t,bt^{1-2n},bs^{-1}t^{-n})\cdot
W_\mu(zs;q,t,as^{-2}t^{-2n},bs^{-1}t^{-n})
\end{multline}
where $z\in\mathbb{C}^n$ and $W_\lambda$ denotes the symmetric Macdonald function~\cite{CoskunG1} that is defined as follows. 

Let $\lambda=(\lambda_1, \ldots, 
\lambda_n)$ and $\mu=(\mu_1, \ldots, \mu_n)$ be partitions of at most
$n$ parts for a positive integer $n$ such that the
skew partition $\lambda/\mu$ is a horizontal strip; i.e. $\lambda_1
\geq \mu_1 \geq\lambda_2 \geq \mu_2 \geq \ldots \lambda_n \geq
\mu_n \geq \lambda_{n+1} = \mu_{n+1} = 0$. Setting
\begin{multline}
\label{definitionHfactor}
H_{\lambda/\mu}(q,p,t,b) \\
:= \prod_{1\leq i < j\leq
n}\left\{\dfrac{(q^{\mu_i-\mu_{j-1}}t^{j-i})_{\mu_{j-1}-\lambda_j}
(q^{\lambda_i+\lambda_j}t^{3-j-i}b)_{\mu_{j-1}-\lambda_j}}
{(q^{\mu_i-\mu_{j-1}+1}t^{j-i-1})_{\mu_{j-1}-\lambda_j}(q^{\lambda_i
    +\lambda_j+1}t^{2-j-i}b)_{\mu_{j-1}-\lambda_j}}\right.\\
\left.\cdot 
\dfrac{(q^{\lambda_i-\mu_{j-1}+1}t^{j-i-1})_{\mu_{j-1}-\lambda_j}}
{(q^{\lambda_i-\mu_{j-1}}t^{j-i})_{\mu_{j-1}-\lambda_j}}\right\}\cdot\prod_{1\leq
i <(j-1)\leq n}
\dfrac{(q^{\mu_i+\lambda_j+1}t^{1-j-i}b)_{\mu_{j-1}-\lambda_j}}
{(q^{\mu_i+\lambda_j}t^{2-j-i}b)_{\mu_{j-1}-\lambda_j}}
\end{multline}
we define
\begin{multline}
\label{definitionSkewW}
W_{\lambda/\mu}(x; q,p,t,a,b)
:= H_{\lambda/\mu}(q,p,t,b)\cdot\dfrac{(x^{-1}, ax)_\lambda
  (qbx/t, qb/(axt))_\mu}
{(x^{-1}, ax)_\mu (qbx, qb/(ax))_\lambda}\\
\cdot\prod_{i=1}^n\left\{\dfrac{\theta(bt^{1-2i}q^{2\mu_i})}{\theta(bt^{1-2i})}
  \dfrac{(bt^{1-2i})_{\mu_i+\lambda_{i+1}}}
{(bqt^{-2i})_{\mu_i+\lambda_{i+1}}}\cdot
t^{i(\mu_i-\lambda_{i+1})}\right\}
\end{multline}
where $q,p,t,x,a,b\in\mathbb{C}$. Note that $W_{\lambda/\mu}(x; q,p, t, a,b)$ vanishes unless $\lambda/\mu$
is a horizontal strip. The symmetric function
$W_{\lambda/\mu}(y, z_1, \ldots, 
z_\ell; q,p,t,a,b)$ is extended to $\ell+1$ variables $y, z_1, \ldots, z_\ell
\in\mathbb{C}$  
through the following recursion formula
\begin{multline}
\label{eqWrecurrence}
W_{\lambda/\mu}(y,z_1,z_2,\ldots,z_\ell;q, p, t, a, b) \\
= \sum_{\nu\prec \lambda} W_{\lambda/\nu}(yt^{-\ell};q, p, t, at^{2\ell},
bt^\ell) \, W_{\nu/\mu}(z_1,\ldots, z_\ell;q, p, t, a, b).
\end{multline}
These functions generalize Macdonald polynomials and interpolation Macdonald polynomials~\cite{Macdonald3, Okounkov1} and are closely related to $BC_n$ abelian functions~\cite{Rains1}. 

\section{Multilateral Basic Series Identities}
I would like to present our bilateralization argument first in one dimensional case to help make the multiple multilateral analogues easier to read. 
The classical \Js, that is $_8\phi_7$ summation formula, for example, may be written in the form 
\begin{multline}
\sum_{k=0}^n \dfrac{(1-bq^{2k})}{(1-b)} \dfrac{(b, q^{-n})_k}{(q, bq^{1+n})_k} 
\dfrac{(\sigma, \rho, \gamma, b^2 q^{1+n}/\sigma \rho \gamma )_k  } 
{(qb/\sigma,qb/\rho, qb/\gamma, \sigma \rho \gamma q^{-n}/b)_{k} } \,q^k \\
= \dfrac{(qb, qb/\sigma\rho, qb/\sigma\gamma, qb/\rho\gamma)_n }{(qb/\sigma, qb/\rho, qb/\gamma,  qb/\sigma\rho\gamma)_n } 
\end{multline}
By using the definition~(\ref{qPochSymbol}) of \qPs\ and the identity 
\begin{equation}
  (a)_k = \dfrac{ (-a)^k q^{\binom{k}{2} } } { (q/a)_{-k}}
\end{equation}
we may flip factors and write the summand in the \lhs\ as
\begin{multline}
\dfrac{q^{-z^2}}{ (q^{-2z} )_{\infty}  (q^{2z})_\infty }
\dfrac{ (\sigma, \rho, \sigma q^{-2z}, \rho q^{-2z})_{\infty} }
{(q^{1-2z}, q^{1+n}, q, q^{1+n+2z})_{\infty} } \\
\cdot (\gamma, q^{4z+1+n}/\sigma \rho \gamma , \gamma q^{-2z}, q^{1+n+2z}/\sigma \rho \gamma)_{\infty}  \hspace*{1.3in} \\
\cdot \dfrac{(q^{1-z-(z+k)}, q^{1+n+z-(z+k)}, q^{1-z+(z+k)}, q^{1+n+z+(z+k)})_{\infty} }
{ (\sigma q^{-z+(k+z) }, \rho q^{-z+(k+z) }, \sigma q^{-z-(k+z)}, \rho q^{-z-(k+z)})_{\infty} } \\
\cdot \dfrac{ q^{(z+k)^2 } (q^{-2(z+k) })_{\infty}  (q^{2(z+k)})_\infty } 
{(\gamma q^{-z+(k+z)},  q^{3z+1+n+(z+k)}/\sigma \rho \gamma , \gamma q^{-z-(k+z)}, q^{3z+1+n-(k+z)}/\sigma \rho \gamma)_{\infty} } 
\end{multline}
where we also set $b=q^{2z}$ for some $z\in\mathbb{C}$. It is clear that the summand is invariant under the maps $(z+k)\leftrightarrow w(z+k)$ for all $w$ in the hyperoctahedral group of rank 1, namely $\mathbb{Z}_2$. It is clear that these maps generate full weight lattice $\mathbb{Z}$ for the root system $C_1$ if $z=m/2$ as illustrated in~\cite{Coskun2} for \RSi. However, we will only consider the case when $m=\delta\in\{0,1\}$, that is $b=q^\delta$. 

Recall also that the Macdonald polynomial identity~\cite{Macdonald4} for the root system $C_1$ of rank 1 may be written as
\begin{equation}
1= \dfrac{1}{1-x^2} + \dfrac{1}{1-x^{-2}}
\end{equation}
By setting $x=q^{(z+k)}$ and multiplying the sum on the \lhs\ by the polynomial identity and simplifying, we get
\begin{multline}
\sum_{k=0}^n \dfrac{(1-q^{\delta+2k})}{(1-b)} \dfrac{(q^\delta, q^{-n})_k}{(q, q^{\delta+1+n})_k} 
\dfrac{(\sigma, \rho, \gamma, q^{2\delta+1+n}/\sigma \rho \gamma )_k  } 
{(q^{1+\delta}/\sigma,q^{1+\delta}/\rho, q^{1+\delta}/\gamma, \sigma \rho \gamma q^{-\delta-n})_{k} } \,q^k \\
= \sum_{k=-n-\delta}^n f(\delta) \dfrac{(q^{-n})_k}{(q^{\delta+1+n})_k} 
\dfrac{(\sigma, \rho, \gamma, q^{2\delta+1+n}/\sigma \rho \gamma )_k  } 
{(q^{1+\delta}/\sigma,q^{1+\delta}/\rho, q^{1+\delta}/\gamma, \sigma \rho \gamma q^{-\delta-n})_{k} } \,q^k \\
=\dfrac{(q^{1+\delta}, q^{1+\delta}/\sigma\rho, q^{1+\delta}/\sigma\gamma, q^{1+\delta}/\rho\gamma)_n }
 {(q^{1+\delta}/\sigma, q^{1+\delta}/\rho, q^{1+\delta}/\gamma,  q^{1+\delta}/\sigma\rho\gamma)_n } 
 \end{multline}
where
\begin{equation}
f(\delta)= \left\{
\begin{matrix}
1 & \mathrm{when}\; \delta=0 \\
1/(1-q^\delta) & \mathrm{when}\; \delta=1
\end{matrix} \right.
\end{equation}
Now we send $n\rightarrow\infty$ applying the dominated convergence theorem for infinite series~\cite{Coskun3} to get 
\begin{multline}
\sum_{k=-\infty}^\infty f(\delta) \, \left(\! \dfrac{q^{(\delta+1)}}{\sigma \rho \gamma } \!\right)^k \!
\dfrac{(\sigma, \rho, \gamma )_k  } 
{(q^{1+\delta}/\sigma,q^{1+\delta}/\rho, q^{1+\delta}/\gamma)_k }
\\ = \dfrac{(q^{1+\delta}, q^{1+\delta}/\sigma\rho, q^{1+\delta}/\sigma\gamma, q^{1+\delta}/\rho\gamma)_\infty }
 {(q^{1+\delta}/\sigma, q^{1+\delta}/\rho, q^{1+\delta}/\gamma,  q^{1+\delta}/\sigma\rho\gamma)_\infty } 
\end{multline}
Here we also used the limit rule 
\begin{equation}
\lim_{a\rightarrow\, 0} a^{k} (x/a)_k = (-1)^k x^k q^{\binom{k}{2}}  
\end{equation}
Now, setting $\delta=0$ gives Bailey's $_3\psi_3$ bilateral summation formula. The $\delta=1$ case appears to be a new bilateral sum. 

We now give a multilateral analogue of Bailey's $_3\psi_3$ bilateral summation formula. 
Recall~\cite{CoskunG1} that when $z=xt^{\delta}$ for some $x\in\mathbb{C}$ and $t^{\delta}=(t^{n-1},t^{n-2},\ldots, t, 1)$, the multiple \Js~(\ref{WJackson}) may be written as
\begin{multline}
\label{simplifiedJs2}
  \dfrac{(sx^{-1}, as x)_\lambda}{(qb x, q b /ax )_\lambda}
 =\sum_{\mu \subseteq \lambda}\,  q^{|\mu|} t^{2n(\mu)} \dfrac{(s, a
   s)_\lambda} {(qb, qb /a)_\lambda} \dfrac{ (bt^{1-n}, qb/a s)_\mu }{
   (qt^{n-1}, as)_\mu} \\ \cdot\prod_{i=1}^n \left\{\dfrac{(1-b
   t^{2-2i}q^{2\mu_i})} {(1-b t^{2-2i})} \right\}
 \prod_{1\leq i<j \leq 
 n} \left\{\dfrac{ (qt^{j-i})_{\mu_i-\mu_j} (bt^{3-i-j})_{\mu_i+\mu_j}}
 {(qt^{j-i-1})_{\mu_i-\mu_j} (b t^{2-i-j})_{\mu_i+\mu_j} }
 \right\}  \\ \cdot W_\mu(q^\lambda t^{\delta(n)}; q, t, b st^{2-2n},
 bt^{1-n}) \dfrac{(x^{-1}, ax)_\mu} {(qbx, q b /ax)_\mu} 
\end{multline}
It was also shown~\cite{Coskun1} that the summand which includes the $W_\mu$ function 
is invariant under the \hgr\ action of permutations an sign changes for arbitrary partitions $\lambda$. More precisely, it was shown that under the specialization $t=q^k$ and $b=q^{2z_i+2k(i-1)}$ where $z_i\in\mathbb{C}$ and $k\geq 0$ is a non--negative integer, the summand is invariant under the action $(\mu_i+z_i)\leftrightarrow w(\mu_i+z_i)$ for all elements $w\in W$, the \hgr\ or rank $n$. It was further verified that this action generates the full weight lattice $\mathbb{Z}^n$ only if $z_i=m/2+k(n-i)$ for some non--negative integers $m, k\geq 0$. 

The proof of the invariance follows from the duality formula and the flip identity~\cite{Coskun2} for $W_{\lambda}$ functions. The symmetry under permutations are given by duality formula
\begin{multline}
\label{eq:duality}
W_{\lambda}\left(k^{-1}q^\nu t^\delta;q,t,k^2a,kb\right)
\cdot\dfrac{(qbt^{n-1})_\lambda (qb/a)_\lambda}{(k)_\lambda
  (kat^{n-1})_\lambda} \\
\cdot \prod_{1\leq i < j\leq n} \left\{ \dfrac{(t^{j-i})_{\lambda_i
    -\lambda_j}(qa^{\prime}t^{2n-i-j-1})_{\lambda_i+\lambda_j}}
{(t^{j-i+1})_{\lambda_i
-\lambda_j}(qa^{\prime}t^{2n-i-j})_{\lambda_i+\lambda_j}} \right\}\\
=W_{\nu}\left(h^{-1}q^\lambda t^\delta;q, t,h^2a^{\prime},hb\right)
\cdot \dfrac{(qbt^{n-1})_\nu (qb/a^{\prime})_\nu}{(h)_\nu
  (ha^{\prime}t^{n-1})_\nu} \\
\cdot \prod_{1\leq i < j\leq n} \left\{ \dfrac{(t^{j-i})_{\nu_i
      -\nu_j}(qat^{2n-i-j-1})_{\nu_i+\nu_j}} {(t^{j-i+1})_{\nu_i
-\nu_j}(qat^{2n-i-j})_{\nu_i+\nu_j}} \right\}
\end{multline}
where $k=a^{\prime}t^{n-1}/b$ and $h = at^{n-1}/b$. The invariance under sign changes follows from the flip identity
\begin{multline}
a^{|\lambda|} b^{-|\lambda|} q^{-|\lambda|}
t^{-n(\lambda)+(n-1)|\lambda|} W_{\lambda}(x_1^{-1}, \ldots, x_n^{-1}; q^{-1}, p, t^{-1}, a^{-1},
b^{-1})\\
= a^{-|\lambda|} b^{|\lambda|} q^{|\lambda|}
t^{n(\lambda)-(n-1)|\lambda|} W_{\lambda}(x_1, \ldots, x_n; q, p, t,
a, b)
\end{multline}
The invariance of other factors follows immediately form the definition of the \qPs.

We will give the multiple $_3\psi_3$ summation using the specialization $m=\delta\in\{0,1\}$ as in the classical one dimensional case, and for $k=1$ or $t=q$. In other words, we let $t\rightarrow q$ and $b\rightarrow q^{\delta+2(n-1)}$ and write the identity above in the form
\begin{multline}
  \dfrac{(sx^{-1}, as x)_\lambda}{(s, a s)_\lambda} \dfrac {(q^{\delta+2n-1}, q^{\delta+2n-1}/a)_\lambda} {(q^{\delta+2n-1} x, q^{\delta+2n-1} /ax )_\lambda}\\
 =\sum_{\mu \in\mathbb{Z}^n}\, f(\delta)\, q^{|\mu|+2n(\mu)}  \dfrac{ ( q^{\delta+2n-1}/a s)_\mu }{
   ( as)_\mu}  \hspace*{0.9in}  \\ 
   \cdot \prod_{1\leq i<j \leq 
 n} \left\{\dfrac{ (q^{j-i+1})_{\mu_i-\mu_j} (q^{\delta+2n+1-i-j})_{\mu_i+\mu_j}}
 {(q^{j-i})_{\mu_i-\mu_j} (q^{\delta+2n-i-j})_{\mu_i+\mu_j} }
 \right\}  \hspace*{0.5in} \\
 \cdot W_\mu(q^{\lambda+\delta(n)}; q, q, sq^{\delta},
 q^{\delta+(n-1)} ) \dfrac{(x^{-1}, ax)_\mu} 
 {(q^{\delta+2n-1}x, q^{\delta+2n-1} /ax)_\mu} 
\end{multline}
where
\begin{equation}
f(\delta):= 
\dfrac{1}{2^n} \prod_{i=1}^{n-1} \dfrac{ 1  }
   { (1+q^{n-i} )  }, \quad \mathrm{if\;} \delta=0 
\end{equation}
and
\begin{equation}
f(\delta):= 
\dfrac{1}{2^n} \prod_{i=1}^{n} \dfrac{ 1  }
   { (1-q^{1+2n-2i} )  }, \quad \mathrm{if\;} \delta=1 
\end{equation}
Note also that although the series is written over $\mathbb{Z}^n$, it actually terminates from above by $\lambda$ and from below by 
$(-\lambda_i - 2n- 2i+\delta) $.

The analogue of Weyl degree formula~\cite{Coskun1} for $W_\mu$ functions
implies that
\begin{multline}
W_{\mu}(q^{N+\delta(n)};q,q, sq^{\delta},q^{\delta+(n-1)}) \\
=\dfrac{(q^{-N},  sq^{\delta+N+n-1})_\mu}{(q^{N+\delta+2n-1}, q^{-N+n}/ s)_\mu} 
\! \prod_{1\leq i < j\leq n} \dfrac{(q^{j-i+1})_{\mu_i
-\mu_j}(q^{\delta+2n-i-j+1})_{\mu_i+\mu_j}}
{(q^{j-i})_{\mu_i -\mu_j}(q^{\delta+2n-i-j})_{\mu_i+\mu_j}}
\end{multline}
Therefore, by setting $\lambda=N^n=(N,N,\ldots,N)$ and sending $N\rightarrow \infty$ we get
\begin{multline}
  \dfrac{(sx^{-1}, as x)_{\infty^n}}{(s, a s)_{\infty^n}} \dfrac {(q^{\delta+2n-1}, q^{\delta+2n-1}/a)_{\infty^n}} {(q^{\delta+2n-1} x, q^{\delta+2n-1} /ax )_{\infty^n}} \dfrac{1}{2^n f(\delta) } \\
 =\sum_{\mu \in\mathbb{Z}^n}\, q^{(1-n)|\mu|+2n(\mu)} s^{|\mu|} \, 
 \dfrac{ (q^{\delta+2n-1}/a s)_\mu }{
   (as)_\mu}  \dfrac{(x^{-1}, ax)_\mu} 
 {(q^{\delta+2n-1}x, q^{\delta+2n-1} /ax)_\mu}  \\ 
\cdot \!\! \prod_{1\leq i<j \leq n} \left\{
 \dfrac{ (1-q^{j-i+\mu_i -\mu_j} ) } {(1-q^{j-i})^2  }
 \dfrac{ (1-q^{\delta+2n-i-j+\mu_i+\mu_j}) } 
 { (1-q^{\delta+2n-i-j } )^2 } 
 \right\} \hspace*{0.3in}
\end{multline}
This is the multilateral analogue of Bailey's bilateral $_3\psi_3$ summation formula as desired. 

\section{Conclusion}
The multilateralization technique applied here can be used to prove other multilateral series when the invariance property is satisfied. We will explore similar identities in other publications.

\end{document}